# Nonstandard Hopf bifurcation in switched systems

Xiao-Song Yang, Songmei Huan

Department of Mathematics,
Huazhong University of Science and Technology.
Wuhan, 430074, China

*Abstract* This paper presents an analysis on nonstandard generalized Hopf bifurcation in a class of switched systems where the lost of stability of linearized systems is not due to the crossing of their complex conjugate eigenvalues but relevant to the switching laws between the subsystems. The mechanism of Hopf-like bifurcation dealt with in this paper is remarkably different from the mechanism of the classical Hopf bifurcation and is also not the same as the generalized Hopf bifurcation studied in the literature.

## 1 Introduction

Piecewise smooth dynamical systems are an important class of ordinary differential equations that arise in scientific problems and engineering application. One of the most interesting class of piece-wise smooth systems are the switched systems. By a switched system we mean a family of continuous time dynamical systems and a rule that determines at each time which dynamical system is responsible for the time evolution. Up to now, much attention has been paid to stability, chaos and control of switched systems [2], [3], [10], [11]. Here we are more concerned with various dynamics of switched systems from dynamics point of view. More precisely we are interested in periodic trajectories and Hopf-like bifurcation in switched systems in this paper. Such topic can be conventionally categorized as Hopf-like bifurcation in piecewise smooth systems.

Many bifurcation properties caused by a discontinuity of dynamical system have so far been widely studied in the literature. These new bifurcation phenomena are unique to piece-wise smooth systems, and are called discontinuity-induced bifurcation (DIB) [1]. Also the term "C-bifurcation" was referred to these discontinuity-induced bifurcation phenomena by Feigin [6] when he studied the doubling of the oscillation period in piecewise continuous systems.

There are various bifurcation phenomena induced by discontinuity. Many researchers have paid attention to the study and classification of different types of DIB in piecewise smooth dynamical systems. Border-collision bifurcation of fixed points in maps refers to the phenomena when a family of fixed points transversely pass across the line of discontinuity as the parameter varies. Grazing bifurcation of limit cycles concerns the corresponding properties when a periodic orbit tangentially intersects the line of discontinuity, which gives rise to complex transitions in piecewise smooth systems. Sometimes, this bifurcation can be analyzed by studying the border-collision bifurcation of a suitably defined Poincaré map near the periodic orbit. Sliding bifurcation appears when part of a periodic orbit coincides with the line of discontinuity, which has vast application backgrounds. For a comprehensive review of all the above phenomena and their classification, the reader is referred to [1],[8],[9] and references there in.

Another type of discontinuity-induced bifurcation concerns the appearance or disappearance of a periodic orbit that is related to Hopf bifurcation or generalized Hopf bifurcation as studied in [4], [5] where a case of Hopf bifurcation in which the periodic orbits are generated from infinity was considered. The paper [8] considered the so-called nonsmooth Hopf bifurcation that an equilibrium hits the line of discontinuity, then becomes unstable and creates periodic oscillations. In papers [7], [13], [14], the authors studied the generalized Hopf bifurcation for the case that the equilibria always stay in the smooth line of discontinuity and the bifurcating periodic orbit transversely crosses the line of discontinuity at least twice. In particular, a generalized Hopf bifurcation emanated from a corner with respect to several discontinuity boundaries was fully discussed [14].

To give a flavor of the DIB-type Hopf bifurcation let us first recall the classical Hopf bifurcation, which is one of the most fundamental mechanisms in bifurcation theory of smooth dynamical systems. Geometrically speaking, Hopf bifurcation means a creation of periodic trajectories in a smooth dynamical systems wherein the linearized system undergoes a change of the phase portrait from a stable node to an unstable node through a center that corresponds to a simple crossing of precisely one pair of complex conjugate eigenvalues of the linearized system through the imaginary axis.

Since the properties of the vector field defined on the line of discontinuity can result in very complicated dynamics and rich bifurcation phenomena in context of piecewise smooth systems, the well known property that a pair of complex eigenvalues of the linearized part go across the imaginary axis is no longer a necessary condition for studying the so called generalized Hopf bifurcation studied in the literature, but other factors in piecewise smooth system have to be taken into account, which reflect the complexity of the bifurcation mechanism for piece-wise smooth systems compared to smooth systems.

In this paper we will present a nonstandard Hopf bifurcation in a class of switched systems that are also called Filippov systems since they are only piece-wise continuous. Here for convenience the terms "nonstandard Hopf bifurcation" refer to the case that the creation of periodic trajectory is not due to the crossing of complex conjugate eigenvalues of the linearized systems through the imaginary axis but other mechanism. The nonstandard Hopf bifurcation considered in this paper is different from the mechanism of the classical Hopf bifurcation but similar to the so-called generalized Hopf bifurcation as well in the literature [14].

The switched systems we consider can be described as follows:

$$\dot{x} = f_{p(t)}(x), \tag{1}$$

$$p(t) = \zeta(x(t), p(t^-)), \tag{2}$$

where $x = (x_1, x_2, ..., x_n)^T \in R^n$ is the state vector, $p(t) \in Q = \{1, 2, ..., m\}$ is the switching index (depending on time), the function $f_i(x): R^n - R^n$ is continuous for every $i \in Q$, and $\zeta(x(t), p(t^-))$ is a transition function depending on the state vector $x(t)$ over the discrete integer set $Q$, in which it is assumed that $p(t)$ is a right-hand continuous function, so $p(t^-)$ is the right-hand limiting value at $t = t^-$.

Systems (1) with the transition law (2) is usually called an autonomous hybrid system [2]. In this paper we discuss switched planar system,

$$\dot{x} = f_{p(t)}(x), \tag{3}$$

$$p(t) = \zeta(x(t), p(t^-)), \tag{4}$$

where $x = (x_1, x_2)^T \in R^2$, $p(t) \in Q = \{1,2,3,4\}\}$. Define four regions

$$S_1 = \{(x_1, x_2) \in R^2 : 0 < x_1 < \infty, x_2 \geq 0\}$$

$$S_2 = \{(x_1, x_2) \in R^2 : -\infty < x_1 \leq 0, x_2 > 0\}$$

$$S_3 = \{(x_1, x_2) \in R^2 : -\infty < x_1 < 0, x_2 \leq 0\}$$

$$S_1 = \{(x_1, x_2) \in R^2 : 0 \leq x_1 < \infty, x_2 < 0\}$$

Clearly

$$R^2 = S_1 \cup S_2 \cup S_3 \cup S_4 \cup \{(0,0)\}$$

The transition function $p(t)$ can be defined as follows in the same way as [2]

$$p(t) = \zeta(x(t), p(t^-)) = i \quad x(t) \in S_i, \quad i = 1,2,3,4 \tag{5}$$

Every half coordinate line (for example the positive $x_1$-half line $\{0 < x_1 < \infty, x_2 = 0\}$) is a switching manifold, where the transition law or switching law (5) takes effect.

## 2 Bifurcation of switched linear systems

In order to study nonstandard Hopf bifurcation of switched nonlinear systems, it is necessary to first investigate the bifurcation problem in the case of switched linear systems

Consider the following switched linear systems

$$\dot{x} = A_{p(t)} \cdot x \tag{6}$$

$$p(t) = \zeta(x(t), p(t^-)), \tag{7}$$

with $A_1 = A_3 = A$ and $A_2 = A_4 = B$, where

$$A = \begin{bmatrix} -a & b \\ -c & -a \end{bmatrix}, \qquad B = \begin{bmatrix} -a & c \\ -b & -a \end{bmatrix} \qquad (8)$$

With the transition function defined in (5), the switched system (6) and (7) can be described equivalently as piecewise linear system

$$\dot{x} = \begin{cases} Ax, & x \in S_1 \cup S_3 \\ Bx, & x \in S_2 \cup S_4 \end{cases} \qquad (9)$$

Throughout this paper we assume that $a > 0, b > 0$ and $c > 0$. It is easy to see that every eigenvalue of the matrix $A$ is a complex number with negative real part, the same is true for $B$. It is easy to see that every subsystem of (6) is asymptotically stable with matrices defined by (8) under the assumption above. In this section we will be studying the behavior of switched system (9).

For this purpose let us first consider the following linear system

$$\dot{x} = Ax$$

The eigeinvalues of $A$ are $\lambda_A = -a \pm i\sqrt{bc}$, and the solution with a given initial point $x_0 = (x_{10}\ x_{20})^T$ has the following form

$$\begin{cases} x_1(t) = \sqrt{x_{10}^2 + \dfrac{b}{c} \cdot x_{20}^2} \cdot e^{-at} \cdot \sin\left(\sqrt{bc} \cdot t + \arctan(\sqrt{\dfrac{c}{b}} \cdot \dfrac{x_{10}}{x_{20}})\right), \\ x_2(t) = \sqrt{\dfrac{c}{b} \cdot x_{10}^2 + x_{20}^2} \cdot e^{-at} \cdot \cos\left(\sqrt{bc} \cdot t + \arctan(\sqrt{\dfrac{c}{b}} \cdot \dfrac{x_{10}}{x_{20}})\right). \end{cases}$$

It can be seen from this expression that the corresponding trajectory $\Gamma_A$ approaches clock-wisely the origin as $t \to +\infty$. Furthermore, it is easy to see that if $b/c$ is large enough, then the portrait of $\Gamma_A$ is narrow in $x_2$ direction, as shown in Fig. 1

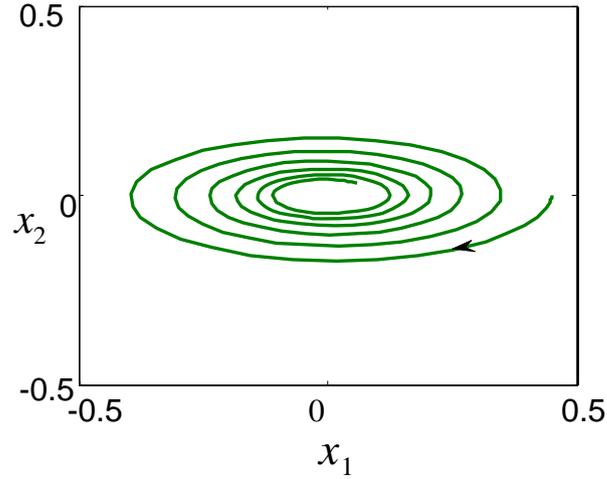

Fig. 1 The portrait of $\dot{x} = Ax$ when $b = 6, c = 1$ and $a = 0.1$

Similarly, for the system
$$\dot{x} = Bx$$

The eigeinvalues of $B$ are $\lambda_B = -a \pm i\sqrt{bc}$, and the solution with a given initial point $x_0 = (x_{10}\ x_{20})^T$ has the following form

$$\begin{cases} x_1(t) = \sqrt{x_{10}^2 + \frac{c}{b} \cdot x_{20}^2} \cdot e^{-at} \cdot \sin\left(\sqrt{bc} \cdot t + \arctan(\sqrt{\frac{b}{c}} \cdot \frac{x_{10}}{x_{20}})\right), \\ x_2(t) = \sqrt{\frac{b}{c} \cdot x_{10}^2 + x_{20}^2} \cdot e^{-at} \cdot \cos\left(\sqrt{bc} \cdot t + \arctan(\sqrt{\frac{b}{c}} \cdot \frac{x_{10}}{x_{20}})\right). \end{cases}$$

The corresponding trajectory $\Gamma_B$ goes clockwisely to the origin as $t \to +\infty$. Furthermore, it is easy to see that if $b/c$ is large enough, then the portrait of $\Gamma_B$ is narrow in $x_1$ direction, as shown in Fig. 2

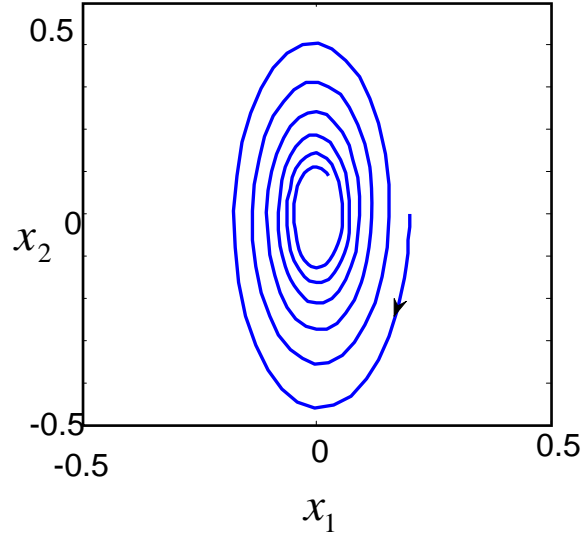

Fig. 2 The portrait of $\dot{x} = Bx$ when $b = 6, c = 1$ and $a = 0.1$

Now we consider the section map $\pi_1$ from the positive semi $x_2$-axis to positive semi $x_1$-axis induced by the flow of $\dot{x} = Ax$. Take point $p = (0, x_2)$ with $x_2 > 0$, the flow $\Gamma_1$ with this initial point is defined by the following solution

$$\begin{cases} x_1(t) = x_2 \sqrt{b/c} \cdot e^{-at} \cdot \sin\left(\sqrt{bc} \cdot t\right) \\ x_2(t) = x_2 \cdot e^{-at} \cdot \cos\left(\sqrt{bc} \cdot t\right) \end{cases}.$$

Thus

$$\pi_1(p) = \pi_1(0, x_2) = (x_2 \sqrt{b/c} \cdot e^{-\frac{a}{\sqrt{bc}} \frac{\pi}{2}}, \quad 0), \quad x_2 > 0$$

Then we consider the section map $\pi_3$ from the negative semi $x_2$-axis to negative semi $x_1$-axis induced by the flow of $\dot{x} = Ax$. Take point $p = (0, x_2)$ with $x_2 < 0$, the flow $\Gamma_3$, i.e., $\Gamma_A$ with this initial point also defines the section map $\pi_3$ as

$$\pi_3(p) = (x_2 \sqrt{b/c} \cdot e^{-\frac{a}{\sqrt{bc}} \frac{\pi}{2}}, \quad 0), \quad x_2 < 0$$

The section map $\pi_2$ from the positive semi $x_1$-axis to negative semi $x_2$-axis and the section map $\pi_4$ from the negative semi $x_1$-axis to positive semi $x_2$-axis induced by the flow of $\dot{x} = Bx$ can be given in the same way as follows,

$$\pi_2(p) = \pi_2(x_1, 0) = (0, \quad -x_1\sqrt{b/c} \cdot e^{-\frac{a}{\sqrt{bc}}\frac{\pi}{2}}), \quad x_1 < 0$$

and

$$\pi_4(p) = \pi_4(x_1, 0) = (0, \quad -x_1\sqrt{b/c} \cdot e^{-\frac{a}{\sqrt{bc}}\frac{\pi}{2}}), \quad x_1 > 0$$

Finally, the Poincaré map $\pi$ from the positive semi $x_1$-axis to itself is a composition of the above section maps and is defined as

$$\pi = \pi_1 \circ \pi_2 \circ \pi_3 \circ \pi_4.$$

Fig.3 is an illustration of the Poincaré map $\pi$.

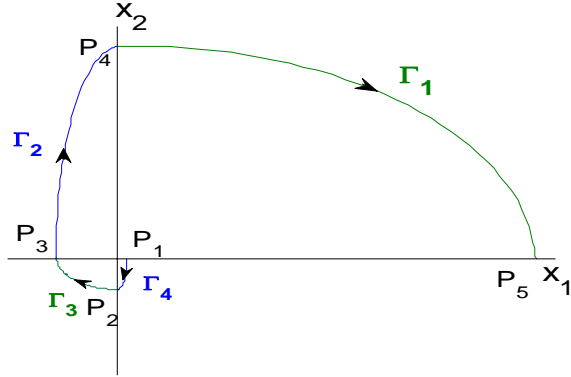

Fig.3   An illustration of the Poincaré map $\pi$.

The brilliant feature of (6), (7) with (8) is that every linear system is asymptotically stable, but the overall switched systems are not stable as shown in Fig.3. A specific example can be found in [2]

Let

$$\Delta = (b/c)^2 \cdot e^{-\frac{a}{\sqrt{bc}} \cdot 2\pi},$$

we have the following theorem.

**Theorem 2.1** For system (6)-(8) with the transition law (5), the following holds:

(i)   There exists a family of periodic solutions if and only if $\Delta = 1$;
(ii) The origin is asymptotically stable if and only if $\Delta < 1$;
(iii) The origin is unstable if and only if $\Delta > 1$.

**Proof** It is easy to see that the Poincaré map $\pi$ from the positive semi $x_1$-axis to itself can be

written as

$$\pi(x_1,0) = (x_1 \cdot \left(\frac{b}{c}\right)^2 \cdot e^{-\frac{a}{\sqrt{bc}} \cdot 2\pi}, 0), \quad x_1 \geq 0.$$

The statements of the theorem can be easily proved by virtue of this formula. □

**Remark 2.1** The expression for the Poincaré map $\pi$ also hold for $x_1 \leq 0$, i.e.,

$$\pi(x_1,0) = (x_1 \cdot \left(\frac{b}{c}\right)^2 \cdot e^{-\frac{a}{\sqrt{bc}} \cdot 2\pi}, 0), x_1 \leq 0$$

## 3 Nonstandard Hopf bifurcation in switched systems

Consider the switched systems

$$\dot{x} = f_{p(t)}(x, \lambda) \qquad (10)$$

where $(x, \lambda) \in R^2 \times I$, $I \subset R$ is an open interval containing zero, with $f_i(x, \lambda)$ being expressed as

$$f_i(x, \lambda) = A_i(\lambda)x + \bar{f}_i(x, \lambda), \quad i = 1,2,3,4,$$

and $\bar{f}_i(x, \lambda)$ is smooth function with $\bar{f}_i(x, \lambda) = o(\|x\|)$, $i = 1,2,3,4$. The transition function $p(t)$ is defined by (5).

In addition, $A_1(\lambda) = A_3(\lambda) = A(\lambda)$ and $A_2(\lambda) = A_4(\lambda) = B(\lambda)$ are defined as follows

$$A(\lambda) = \begin{bmatrix} -a & b(\lambda) \\ -c(\lambda) & -a \end{bmatrix}, \qquad B = \begin{bmatrix} -a & c(\lambda) \\ -b(\lambda) & -a \end{bmatrix}$$

where $a > 0, b(\lambda) > 0$ and $c(\lambda) > 0$.

Let $R^+ = \{(x_1, 0) \in R^2 : x_1 \geq 0\}$ and define a function

$$\Delta(\lambda) = \left(b(\lambda)/c(\lambda)\right)^2 \cdot e^{-\frac{a}{\sqrt{b(\lambda)c(\lambda)}} \cdot 2\pi}.$$

We have the following theorem.

**Theorem 3.1** Suppose that every $\bar{f}_i(x, \lambda)$ is analytic and

$$\Delta(0) = 1, \quad \left.\frac{d\Delta(\lambda)}{d\lambda}\right|_{\lambda=0} \neq 0, \qquad \lambda \in I$$

There is an open interval $\bar{I} \subset R$ with $0 \in \bar{I}$ and a continuous function $\bar{\lambda} : \bar{I} \to I$ satisfying

$\bar{\lambda}(0) = 0$ such that for each $x_1 \in \bar{I}$ there is a period orbit of (10) passing though $(x_1, 0)$ at the parameter $\bar{\lambda}(x_1)$.

***Proof*** It is enough to study the Poincaré map $\pi(x_1, \lambda) : R^+ \times I \to R^+$ derived from (10) in a neighborhood of $(x, \lambda) = (0, 0)$.

First we consider the section map $\pi_1$ from the positive semi $x_2$-axis to positive semi $x_1$-axis induced by $f_1(x, \lambda)$, as argued in [12], it is easy to see from the previous arguments for **Theorem 2.1** that

$$\pi_1(0, x_2) = (x_2\sqrt{b(\lambda)/c(\lambda)} \cdot e^{-\frac{a}{\sqrt{b(\lambda)c(\lambda)}}\frac{\pi}{2}} + h_1(\lambda)x_2^{k_1} + o(x_2^{k_1}),\ 0),\ x_2 > 0.$$

Similarly induced by $f_2(x, \lambda)$ we have

$$\pi_2(x_1, 0) = (0,\ -x_1\sqrt{b(\lambda)/c(\lambda)} \cdot e^{-\frac{a}{\sqrt{b(\lambda)c(\lambda)}}\frac{\pi}{2}} + h_2(\lambda)x_1^{k_2} + o(x_1^{k_2})),\ x_1 < 0.$$

Induced by $f_3(x, \lambda)$ we have

$$\pi_3(0, x_2) = (x_2\sqrt{b(\lambda)/c(\lambda)} \cdot e^{-\frac{a}{\sqrt{b(\lambda)c(\lambda)}}\frac{\pi}{2}} + h_3(\lambda)x_2^{k_3} + o(x_2^{k_3}),\ 0),\ x_2 < 0$$

and induced by $f_4(x, \lambda)$ we have

$$\pi_4(x_1, 0) = (0,\ -x_1\sqrt{b(\lambda)/c(\lambda)} \cdot e^{-\frac{a}{\sqrt{b(\lambda)c(\lambda)}}\frac{\pi}{2}} + h_4(\lambda)x_1^{k_4} + o(x_1^{k_4})),\ x_1 > 0$$

Here every $h_i(\lambda)$ satisfies $h_i(0) \neq 0$ and each $k_i > 1$, for $i = 1, 2, 3, 4$.

Now every fixed parameter $\lambda$ the Poincaré map $\pi(x_1, \lambda) : R^+ \to R^+$ is the composition $\pi = \pi_1 \circ \pi_2 \circ \pi_3 \circ \pi_4$ and can be expressed as

$$\pi(x_1, \lambda) = \Delta(\lambda)x_1 + \delta(\lambda)x_1^k + o(x_1^k),$$

in a neighborhood of $(x, \lambda) = (0, 0)$, where $\delta(0) \neq 0$ and $k = \min\{k_i\} > 1$。

Let $\Pi = \pi(x_1, \lambda) - x_1,\ R^+ \to R^+$,

$$\Pi(x_1, \lambda) = x_1\overline{\Pi}(x_1, \lambda) = x_1[(\Delta(\lambda) - 1) + \delta(\lambda)x_1^{k-1} + o(x_1^{k-1})].$$

Therefore the nontrivial periodic orbits of (10) correspond to the zero points of $\overline{\Pi}(x_1, \lambda)$. Now

$$\overline{\Pi}(0,0) = 0, \quad \left.\frac{\partial \overline{\Pi}(0,\lambda)}{\partial \lambda}\right|_{\lambda=0} = \left.\frac{d\Delta(\lambda)}{d\lambda}\right|_{\lambda=0} \neq 0.$$

By the well known implicit function theorem and **Remark 2.1**, it follows that there is an open interval $\overline{I} \subset R$ with $0 \in \overline{I}$ and a continuous function $\overline{\lambda}: \overline{I} \to I$ satisfying $\overline{\lambda}(0) = 0$ such that for each $x_1 \in \overline{I}$ there is a corresponding parameter $\overline{\lambda}(x_1)$ such that $x_1$ is a fixed point of $\Pi(x_1, \overline{\lambda})$. □

*Remark 3.1* If $\bar{f}_i(x, \lambda)$ is not analytic, then the map $\pi(x_1, \lambda)$ is $C^1$ but not $C^2$ in general as remarked in [13]. In view of the arguments of [7] and [13], we can have the same result.

The above theorem can be further refined. For this purpose we need a simple fact about the implicit function that is easy to prove.

**Lemma** Let $F(x, y)$ be an analytic function satisfying $F(0,0) = 0$ and

$$\left.\frac{\partial F(x,y)}{\partial y}\right|_{(0,0)} \neq 0,$$

$$\left.\frac{\partial^i F(x,y)}{\partial x^i}\right|_{(0,0)} = 0, \; i = 1,..,k-1; \quad \left.\frac{\partial^k F(x,y)}{\partial x^k}\right|_{(0,0)} \neq 0.$$

Then there is a neighborhood of $x = 0$ and an analytic function $y(x)$ defined on this neighborhood such that $F(x, y(x)) = 0$ and $y(x) = \zeta x^k + o(x^k)$, where

$$\zeta = \frac{-1}{k!}\left.\left(\frac{\partial^k F(x,y)}{\partial x^k} \Big/ \frac{\partial F}{\partial y}\right)\right|_{(0,0)}.$$

Now consider the Poincaré map defined in **Theorem 3.1**

$$\pi(x_1, \lambda) = \Delta(\lambda)x_1 + \delta(\lambda)x_1^k + o(x_1^k)$$

In view of this Poincaré map, it is apparently enough to consider the following map

$$\overline{\Pi}(x_1, \lambda) = (\Delta(\lambda) - 1) + \delta(\lambda)x_1^{k-1} + o(x_1^{k-1})$$

in a neighborhood of $(x, \lambda) = (0,0)$, where $\delta(0) \neq 0$ and $k = \min\{k_i\} > 1$.

**Theorem 3.2** In addition to the statements in **Theorem 3.1** we have the following facts:

1) If $\delta(\lambda)\dfrac{d\Delta(\lambda)}{d\lambda}\bigg|_{\lambda=0} < 0$, then there exists $\bar{\lambda} > 0$ and a continuous function

$$x_1(\lambda) : [0, \bar{\lambda}) \to R^+, \text{ with } x_1(0) = 0$$

such that for each $\lambda \in (0, \bar{\lambda})$ there corresponds a $x_1(\lambda) > 0$ such that system (10) has a nontrivial period orbit passing though $(x_1(\lambda), 0)$.

2) If $\delta(\lambda)\dfrac{d\Delta(\lambda)}{d\lambda}\bigg|_{\lambda=0} > 0$, then there exists $\hat{\lambda} > 0$ and a continuous function

$$x_1(\lambda) : (-\hat{\lambda}, 0] \to R^+, \text{ with } x_1(0) = 0.$$

such that for each $\lambda \in (-\hat{\lambda}, 0)$ there corresponds a $x_1(\lambda) > 0$ such that system (10) has a nontrivial period orbit passing though $(x_1(\lambda), 0)$ .

*Proof* In view of the above Lemma, we have

$$\lambda(x_1) = \gamma x_1^{k-1} + o(x_1^{k-1})$$

where

$$\gamma = -(\delta(\lambda)/\Delta'(\lambda))\big|_{\lambda=0} = -\frac{\delta(0)}{\Delta'(0)}.$$

It is easy to prove the above statements by virtue of this formula. □

In the following we present some global results on bifurcation of periodic orbits.
Let $S$ be the skew matrix

$$S = \begin{bmatrix} 0 & -1 \\ 1 & 0 \end{bmatrix}.$$

We have the following result.

**Proposition 3.1** Suppose that the following conditions are satisfied:

1) there exists a positive number $M$ and a positive Lyapunov function $V(x)$ such that

$$\left\langle \frac{\partial V(x)}{\partial x}, f_i(x, \lambda) \right\rangle < 0 \text{ for } \|x\| > M \;;$$

2) $\left|\langle A_i(\lambda)x, Sx \rangle\right| > \left|\langle \bar{f}_i(x, \lambda), Sx \rangle\right|$ for $x \neq 0$, $i = 1, 2, 3, 4$;

3) $\Delta(0) = 1$, $\left.\dfrac{d\Delta(\lambda)}{d\lambda}\right|_{\lambda=0} > 0$, $\quad \lambda \in I$.

Then at $\lambda = 0$ there bifurcates a periodic orbit of (10) circling around the origin.

The proof is easy we just give an outline of the proof.

*Outline of the proof*

Consider an interval $J \subset R^+ = \{(x_1, 0) \in R^2 : x_1 \geq 0\}$ defined as

$$J \subset \{(x_1, 0) \in R^2 : 0 \leq x_1 \leq 2M\}$$

Condition 2) guarantee that every orbit except for the origin of (10) goes around the origin so that the Poincaré map $\pi(x_1, \lambda) : J \to J$ can make sense. Conditions 1) and 3) make sure that for $\lambda > 0$ the Poincaré map $\pi(x_1, \lambda)$ can have a fixed point.

**Corollary** Suppose the following conditions are satisfied:

1) $\langle x, \bar{f}_i(x, \lambda) \rangle < 0$, $\quad \dfrac{\|x\|^2}{\left|\langle x, \bar{f}_i(x, \lambda) \rangle\right|} \to 0$ as $\|x\| \to \infty$, $i = 1,2,3,4$;

2) $\left|\langle A_i(\lambda)x, Sx \rangle\right| > \left|\langle \bar{f}_i(x, \lambda), Sx \rangle\right|$ for $x \neq 0$, $i = 1,2,3,4$;

3) $\Delta(0) = 1$, $\left.\dfrac{d\Delta(\lambda)}{d\lambda}\right|_{\lambda=0} > 0$, $\quad \lambda \in I$.

Then at $\lambda = 0$ there bifurcates a periodic solution circling around the origin.

## 4 An example

In this section we present an example to give an illustration of Theorem 3.1.

Consider the switched systems

$$\dot{x} = f_{p(t)}(x, \lambda), \quad p(t) \in \{1,2,3,4\} \qquad (10)$$

with $f_i(x, \lambda)$ being expressed as

$$f_i(x, \lambda) = A_i(\lambda)x + \bar{f}_i(x, \lambda), \quad i = 1,2,3,4,$$

where

$$\bar{f}_{1,3}(x, \lambda) = \begin{bmatrix} -(x_1^3 + \lambda \cdot x_1 x_2^2) \\ -(\lambda \cdot x_2^3 + x_2 x_1^2) \end{bmatrix},$$

and

$$\bar{f}_{2,4}(x,\lambda) = \begin{bmatrix} -x_1^5 \lambda \\ -x_1^4 x_2 \lambda \end{bmatrix}.$$

The transition function $p(t)$ is defined by (5).

In addition, $A_1(\lambda) = A_3(\lambda) = A(\lambda)$ and $A_2(\lambda) = A_4(\lambda) = B(\lambda)$ are defined as follows

$$A(\lambda) = \begin{bmatrix} -a & b(\lambda) \\ -c(\lambda) & -a \end{bmatrix}, \qquad B = \begin{bmatrix} -a & c(\lambda) \\ -b(\lambda) & -a \end{bmatrix}.$$

Take $a = 2$, $b = e\pi$, $c = \pi/e$, $b(\lambda) = b + \lambda^2 + \lambda$, $c(\lambda) = c + \lambda^2$.

Then

$$\Delta(0) = (b/c)^2 \cdot e^{-\frac{a}{\sqrt{bc}} \cdot 2\pi} = 1,$$

and

$$\Delta'(0) = \left( \frac{b(\lambda)}{[c(\lambda)]^3} \cdot e^{\frac{-2a\pi}{\sqrt{b(\lambda) \cdot c(\lambda)}}} [(\frac{a\pi}{\sqrt{b(\lambda) \cdot c(\lambda)}} + 2) \cdot b'(\lambda) \cdot c(\lambda) + (\frac{a\pi}{\sqrt{b(\lambda) \cdot c(\lambda)}} - 2) \cdot b(\lambda) \cdot c'(\lambda)] \right)\bigg|_{\lambda=0}$$

$$= \frac{4}{e\pi} \cdot b'(0) = \frac{4}{e\pi} > 0.$$

From Theorem **3.1** it can be asserted that there bifurcate periodic orbits for $\lambda > 0$, which are sketched in the following Fig.4.

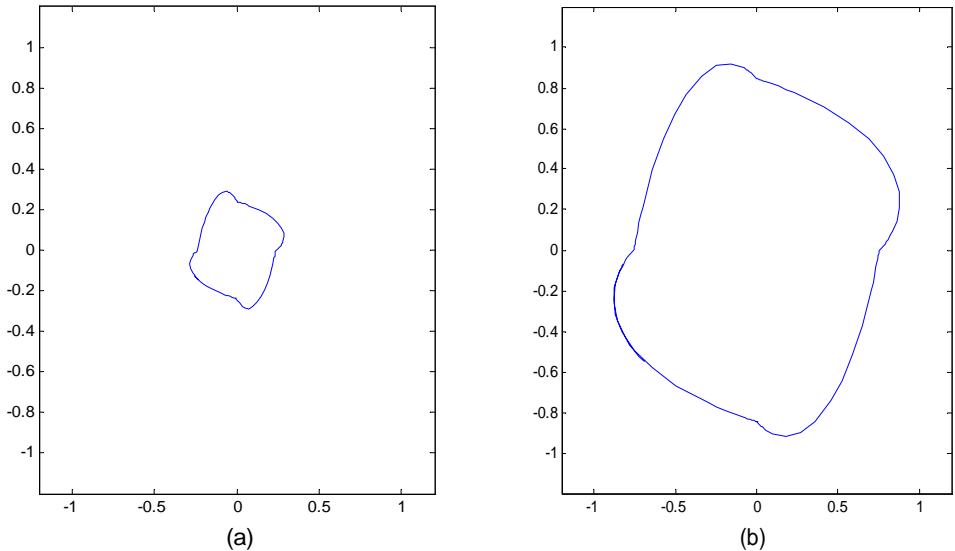

**Fig.4** Creation of periodic orbits for different parameter: (a) $\lambda = 0.1$ (b) $\lambda = 1$

In fact, the creation of a periodic orbit when the parameter varies is a global fact for this example. Since

$$<V_x, A(\lambda)x + \bar{f}_{1,3}(x,\lambda)> = -a \cdot (x_1^2 + x_2^2) + [b(\lambda) - c(\lambda)] \cdot x_1 x_2 - (x_1^4 + \lambda x_2^4 + \lambda x_1^2 x_2^2 + x_1^2 x_2^2)$$

$$<V_x, B(\lambda)x + \bar{f}_{2,4}(x,\lambda)> = -a \cdot (x_1^2 + x_2^2) + [c(\lambda) - b(\lambda)] \cdot x_1 x_2 - x_1^4 (x_2^2 + x_1^2)\lambda \; ;$$

$$|<A(\lambda)x, Sx>| = b(\lambda)x_2^2 + c(\lambda)x_1^2, \qquad <\bar{f}_{1,3}(x,\lambda), Sx> = 0$$

$$|<B(\lambda)x, Sx>| = b(\lambda)x_1^2 + c(\lambda)x_2^2, \qquad <\bar{f}_{2,4}(x,\lambda), Sx> = 0,$$

showing that conditions in **Proposition 3.1** are satisfied globally, therefore the existence of periodic orbit is assured for every $\lambda > 0$.

## Conclusion

In this paper we have presented a nonstandard Hopf bifurcation in a class of switched systems. The lost of stability of the linearized systems and creation of periodic orbits is not due to the crossing of their complex conjugate eigenvalues but closely relevant to the switching laws or transition laws between the subsystems. This is remarkably different from the mechanism of the classical Hopf bifurcation, also it is not the same as the so called generalized Hopf bifurcation.

**Acknowledgements** This work is supported in part by National Natural Science Foundation of China （10972082）.